\documentclass[%
 reprint,
 amsmath,amssymb,
 aps,
]{revtex4-2}
\usepackage{graphicx}% Include figure files
\usepackage{dcolumn}% Align table columns on decimal point
\usepackage{bm}% bold math

\begin{document}

\title{A proof-of-principle experiment on the spontaneous symmetry breaking machine and numerical estimation of its performance on the $K_{2000}$ benchmark problem}

\author{Toshiya Sato$^{1}$}
\email[Toshiya Sato]{:tshy.sato@ntt.com}
\author{Takashi Goh}

\affiliation{NTT Device Technology Laboratories, NTT Corporation, 3-1, Morinosato Wakamiya, Atsugi-shi, Kanagawa, Japan 243-0198}

\date{\today}

\begin{abstract} 
In a previous paper, we proposed a unique physically implemented type simulator for combinatorial optimization problems, called a spontaneous symmetry breaking machine (SSBM). 
In this paper, we first report the results of the experimental verification of SSBM using a small-scale benchmark system, and then describe numerical simulations using the benchmark problems ($K_{2000}$) conducted to confirm its usefulness for large-scale problems.
From 1000 samples with different initial fluctuations, it became clear that SSBM can explore a single extremely stable state.
This unique behavior holds the potential to become a significant advantage over other simulators.
\end{abstract}

\maketitle
%%%%%%%%%%%%%%%%%%%%%%%%%%  body  %%%%%%%%%%%%%%%%%%%%%%%%%%
\section{Introduction}
Symmetry and its (spontaneous) breaking are one of the most significant concepts for understanding various natural phenomena, and it is widely known that they play an essential role in a wide range of areas covered by science \cite{Anderson}. In applied physics, which is at the heart of scientific applications, and in condensed matter physics, which is the main foundation of applied physics, there are many arguments suggesting that spontaneous symmetry breaking (SSB) may play an important role in various phenomena. In those, the discussion is mainly focused on explaining phase transitions (e.g. spontaneous magnetization \cite{Onsager}, insulator-metal phases \cite{Inami}) and changes of state before and after phase transitions, and the main emphasis is on understanding the manifestation of macroscopic properties based on the assumption of the existence of microscopic elements (e.g. electron spin). In dissipative systems, concepts such as dissipative structure and synergetics have been discovered, leading to a broader and deeper understanding of natural phenomena from the perspective of {\it dynamical macroscopic orders} created by symmetry breaking \cite{Prigogine,Haken}. In the field of elementary particle theory, it goes even further and plays an important role in understanding the manifestation of elementary particles (microscopic elements) themselves such as the Higgs boson \cite{Higgs,Nambu}.
  
  In recent years, there has been growing interest in attempts to solve problems in an extremely wide range of fields (logistics, financial portfolios, drug discovery, machine learning, etc.) by treating them as combinatorial optimization problems (COPs).
At the same time, there has been a lot of research focused on Ising problems that can be mapped to COPs in the NP-hard category \cite{Lucas}.
In particular, quantum annealers (a type of quantum computer) \cite{Nishimori,Johnson} and coherent Ising machines (based on laser oscillation phenomena) \cite{Yamamoto} that use physically implemented type simulators to explore the ground state of Ising model systems are attracting attention.

In our previous paper \cite{Toshiy}, we introduced a theoretical framework for understanding a robust causality generated in a unique dissipative system that can be constructed using modern photonic technology. 
Furthermore, based on this dissipative causality, we proposed and demonstrated an experimental system that holds the potential to enable the observation of spontaneous symmetry breaking phenomena and serve as a foundation for its applications.
In further developed discussions, we discovered an unprecedented new phenomenon that can be understood as a more complex SSB phenomenon. We then clarified that this phenomenon can be considered an experimentally verifiable phenomenon possessing duality with the observationally difficult phenomenon \cite{Unruh,Philbin} of creating many-body systems via SSB.
In addition, we proposed a physically implemented simulator (i.e., the spontaneous symmetry breaking machine (SSBM)) to explore good solutions to combinatorial optimization problems that utilize this phenomenon and reported verification results from numerical simulations using a small-scale benchmark problem (${\rm MaxCut^{3}}\!:\!N\!=\!16$).

This paper reports on experimental verification of the SSBM using the  ${\rm MaxCut^{3}} (N\!=\!16)$ problem, \cite{Takata} and, with a view to application to large-scale problems, reports on evaluation by numerical simulation using the $K_{2000}$ problem \cite{Haribara}, which is used as a benchmark problem in various Ising machines \cite{Inagaki,Gotoh01,Gotoh02,KYamamoto}.

\section{SSB machine}
Understanding the dissipative causality underlying the SSB phenomenon observed in SSBM requires identifying dissipative systems with properly defined virtual boundaries \cite{Endnote01}. Fortunately, in SSBM constructed using optical devices like optical fibers and Mach-Zehnder modulators (MZMs), the high optical wave-guiding properties aid in imagining (identifying) suitable boundaries for such dissipative systems.

\begin{figure}[h]
\centering
\includegraphics[width=7cm]{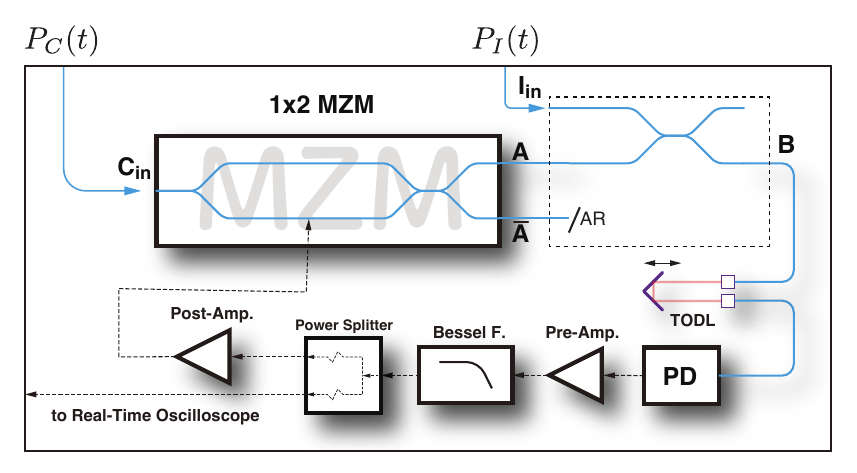}
 \caption{Schematic diagram of the circuit configuration generating dissipative causality. Optical coherent clock pulses $P_{C}(t)$ are supplied to the circuit via the input port ${\bf C}_{\rm in}$, and a 1$\times$2 MZM functions as a partial inflow gate to the system. The electrical pulse originating from the $i$-th optical coherent clock pulse is timed by a tunable optical delay line (TODL) for use in intensity modulation of the ($i\!+\!m$)-th optical coherent clock pulse, and its width is expanded by a Bessel filter to be sufficiently wider than the optical coherent clock pulse width.}
 \end{figure}
Figure 1 is a schematic diagram of the circuit setup that generates dissipative causality, which is the basis of SSBM. Replacing the area enclosed by the dotted line with an optical interference sub-circuit that simulates the pseudo-spin interaction (pSI) to be described later results in a schematic diagram of the SSBM configuration.
Figure 2 is a conceptual diagram of the dynamic system that can be created using this basic circuit, namely the fully dissipative system on which we are focusing.
 \begin{figure}[h]
\centering
\includegraphics[width=7cm]{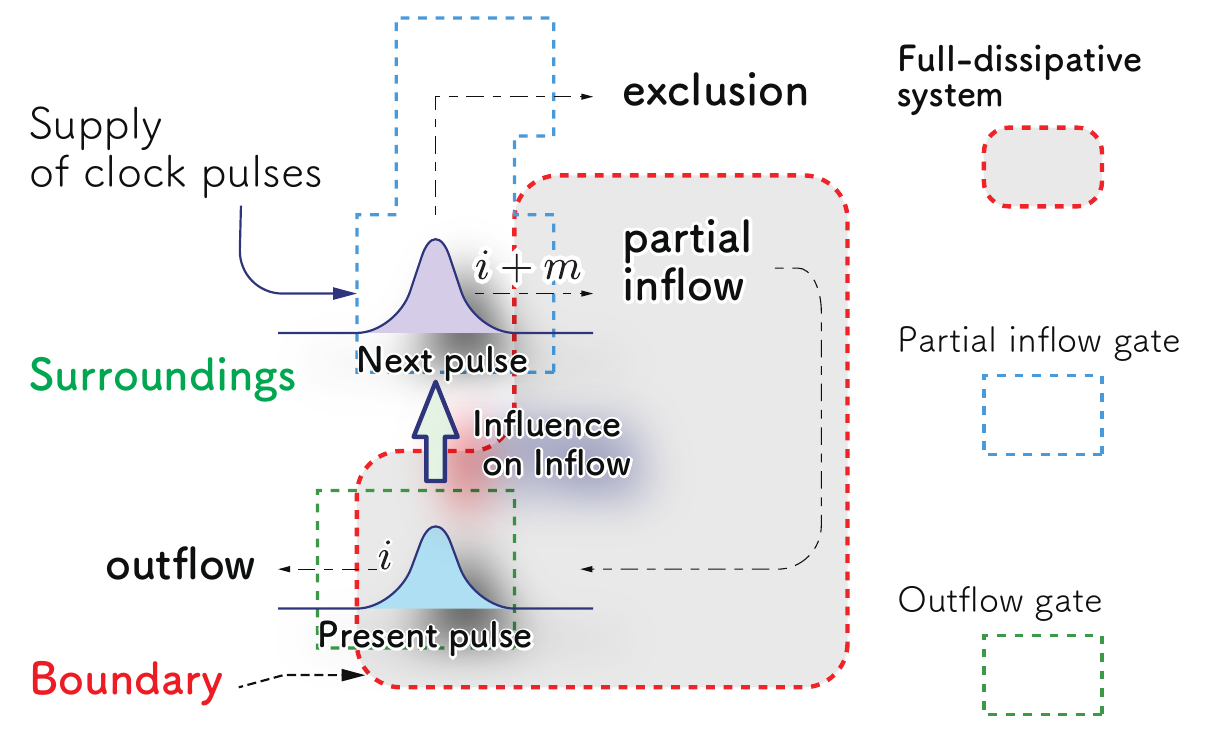}
 \caption{Conceptual diagram of a full-dissipative system. This system possesses appropriate boundaries where the inflow of the ($i\!+\!m$)-th clock pulse is regulated by a partial inflow gate that utilizes the previously inflowed $i$-th clock pulse as an intensity-modulated signal. Simultaneously with the inflow of the ($i\!+\!m$)-th clock pulse, the $i$-th clock pulse completely flows out of the system.}
 \end{figure}
Polarized coherent clock pulses $P_{C}(t)$ with pulse width $\sigma_{\mbox{\rm \scriptsize pw}}$ and repetition interval $\Delta t$ are input to the optical input port ${\bf C}_{{\rm in}}$ of a 1$\times$2 MZM used as a partial inflow gate. The $i$-th output light pulse from output port {\bf A} of the MZM is converted into an electrical pulse by a photodiode after passing through an optical variable delay line to adjust the delay and timing. This electrical pulse is then expanded to a pulse width of $\tau_{\mbox{\rm \scriptsize pw}}$ using a Bessel filter, while its power is compensated by two amplifiers while suppressing excessive noise, and is used for intensity modulation by the MZM for the ($i\!+\!m$)-th optical input pulse.
In the above operation of this circuit, we can identify $m$ independent fully dissipative systems defined by $m$ virtual boundaries where the $i$-th and ($i\!+\!m$)-th clock pulses interact, with the MZM serving as the common inflow gate.
And, when the condition $\sigma_{\mbox{\rm \scriptsize pw}} \ll \tau_{\mbox{\rm \scriptsize pw}} \ll \Delta t$ is satisfied, this fully dissipative system maintains an extremely singular state described by the following iterative equation using the parameters phase difference-induced efficiency $\gamma$, the transmittance $\phi_{i}$ and $\phi_{i+m}$ to port {\bf A} at the input of the clock pulses $i$ and ($i\!+\!m$)-th, and the static phase condition (offset) $\theta_{{\rm B}}$, related to MZM, thus establishing the dissipative causality.
\begin{eqnarray}% equation no. 04
\phi_{i+m} 
& = &
\sin^{2}(\frac{\gamma}{2} \, \phi_{i}  +  \theta_{B} ) 
\end{eqnarray}
Here, we set the remaining degrees of freedom as $\theta_{B}\!=\!0$ and $\gamma\!=\!\pi$. These are the parameters that enable the SSB phenomenon to occur in the SSBM \cite{Toshiy, Endnote04}.
(Unlike examples that realize systems following similar iterative equations by incorporating analog functions into the digital processing flow using MZM characteristics \cite{Bohm}, robust causality as a phenomenon is realized here.)
However, when the optical clock pulse is supplied and the dissipative causality is established, the states $\phi_{i}\!=\!0$ appear and are maintained (becoming hidden states).
Therefore, it is necessary to input an initialization optical pulse train $P_{I}(t)$ 
from the optical input port $\bf I_{{\rm in}}$ with adjusted timing to transition all full-dissipative systems to an unstable fixed point ($\phi_{i}\!=\!1/2$), which is the state before SSB begins (SBSB).
Moreover, we will consider the construction of a many-body-like system by imparting interrelationships between $m$ full-dissipative systems capable of causing the SSB phenomenon described above.
Specifically, we treat the two attractors in the above full-dissipative systems as pseudo spin up and spin down, and introduce pseudo-spin interactions (pSIs) represented by the following equations, which can be implemented by an optical interference circuit.
\begin{eqnarray}% equation no. 02 ~ 7
\phi_{i+m}
& = &
\sin^{2}(\frac{\pi}{2} \,|\sqrt{\phi_{i}} -  {\cal Q}_{i} |^{2} ) \\
{\cal Q}_{i}
& = &  
{\cal Q}_{i}^{\rm F} + {\cal Q}_{i}^{\rm AF} \\
{\cal Q}_{i}^{\rm F}
& = &  
- \sum^{i \ne k} {\cal J}_{i:k}^{\rm F}( \sqrt{\phi_{i}}  -  \sqrt{\phi_{i+k}})  \\
{\cal Q}_{i}^{\rm AF}
& = &  
- \sum^{i \ne k} {\cal J}_{i:k}^{\rm AF}( \sqrt{\phi_{i}}  -  \sqrt{\overline{\phi_{i+k}}})
\end{eqnarray}
\begin{eqnarray}
 {\cal J}_{i:k}^{\rm F} = \left\{ \begin{array}{ll}
  {\cal J}_{i:k} & \mbox{if} \mbox{\footnotesize \rm \,\,\,\,\,\,\,\,\,Ferro.\, type}\\
                                 0 & \mbox{if}\,\,\, \mbox{\footnotesize \rm  Anti-Ferro.\, type} \end{array} \right. \\
 {\cal J}_{i:k}^{\rm AF} = \left\{ \begin{array}{ll}  0 & \mbox{if}\,\,\, \mbox{\footnotesize \rm  \,\,\,\,\,\,\,\,\,Ferro.\, type} \\
                                  {\cal J}_{i:k} & \mbox{if} \,\,\,   \mbox{\footnotesize \rm  Anti-Ferro.\, type}  \end{array} \right. 
\end{eqnarray}
where $k \!\in\! I, -m/2 \!\le\! k \!<\! m/2$.  Note that $\overline{\phi}$ is the transmittance to the port $\overline{\bf A}$ and is also the complement of $\phi$ ($\overline{\phi} = 1 - \phi$).

As is clear from the equations these pSIs are different from the so-called exchange interactions, but they satisfy the behavior consistent with the interpretation of exchange energy (namely, ferromagnetic types stabilize when states are aligned, while antiferromagnetic types stabilize when states are opposite).
This is necessary to deal with the fact that pseudo-spins are continuous variables rather than discrete variables.
Even though pseudo-spins are continuous variables, if we use interactions that follow exchange interactions as they are, we cannot guaranty the expected behavior in regions other than those corresponding to the values assigned to spins.
The introduction of these pSIs yields new and interesting results.
As is clear from Equations (4) and (5), when $m$ full-dissipative systems are simultaneously in the SBSB $(\phi\!=\!1/2)$, the effects of pSIs disappear.
Therefore, it is possible to transition the entire complex system consisting of $m$ full-dissipative systems with pSIs to SBSB through the initialization process described above.
Moreover, this complex system creates a many-body-like system consisting of pseudo-spins, incorporating the effects of pSIs through the SSB phenomenon.
We focus on the phenomenon described above - namely, the spontaneous breaking of symmetry that was temporarily restored by the initialization process, leading to a transition to a stable fixed point (i.e., the SSB phenomenon) - and utilize this as the attractive force inherent in the discrete dynamics itself that drives the transition to two stable fixed points.
We refer to a many-body-like system composed of pseudo-spins as described above as an Ising-model-like system.

\section{Proof-of-principle experiment}

To understand the creation of the Ising-model-like system through the SSB phenomenon, we cannot simply apply the concepts developed in condensed matter physics.
And we do not know of any other specific examples of the phenomenon of many-body-like systems being created by SSB.
However, since consistency with the interpretation of exchange interaction energy was ensured when introducing pSIs, it is expected that a correspondence will exist between the stable states in the Ising model system and those in the Ising-model-like system.
Nevertheless, the correspondence with the Ising model system discussed \cite{Leleu, Bohm02}. regarding continuous-variable Ising machines utilizing bifurcation phenomena through feedback strength control cannot be directly applied (see also Supplementary Materials).
In our previous paper, we conducted numerical simulations to confirm this correspondence and found it promising \cite{Toshiy}.
Here, we present one of the most compelling pieces of evidence to ensure the reliability of this correspondence (duality) by experimentally demonstrating the effectiveness of SSBM as a simulator for COPs.
\begin{figure}[h]
\centering
\includegraphics[width=7cm]{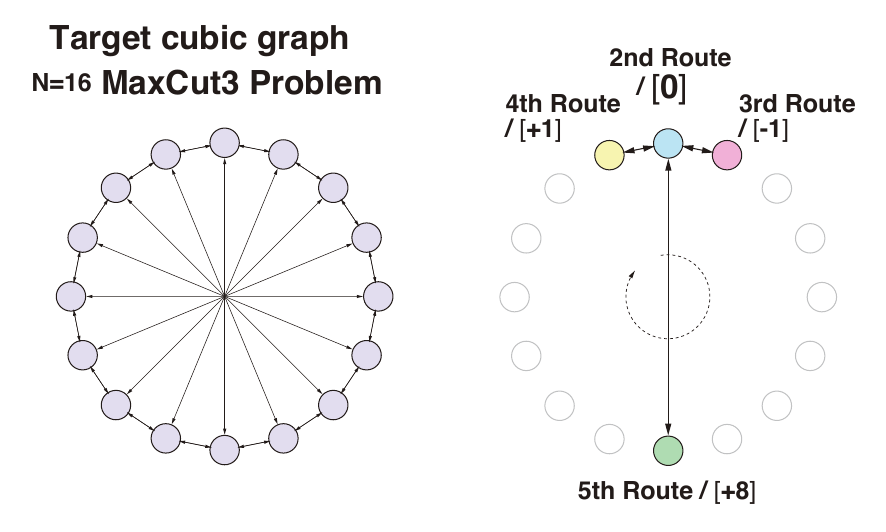}
 \caption{A 3rd-order graph of a target COP ($N\!=\!16$), which is addressed in the proof-of-principle experiment of the SSBM with  ${\cal J}_{i:k}^{\rm AF}\!=\!1/30\, (k\!=\!i\!+\!8, i\!+\!1, i\!-\!1)$ and ${\cal J}_{i:k}^{\rm F}\!=\!0$. This graph represents NP-hard instances, equivalent to MaxCut3 problem.}
 \end{figure}
Figure 3 is a third-order graph $(N\!=\!16)$ representing a COP belonging to the category called MaxCut3 \cite{Takata}, which is addressed in the proof-of-principle experiment.
This COP is an established problem that was used in demonstration experiments of coherent Ising machines, which are one type of physically implemented Ising machine.
And this COP has the advantage that dynamic switching is not required in the optical delay interference circuit (ODIC) \cite{Goh,Endnote02} for the physical implementation of pseudo-spin interactions.
 \begin{figure}[h]
\centering
\includegraphics[width=9cm]{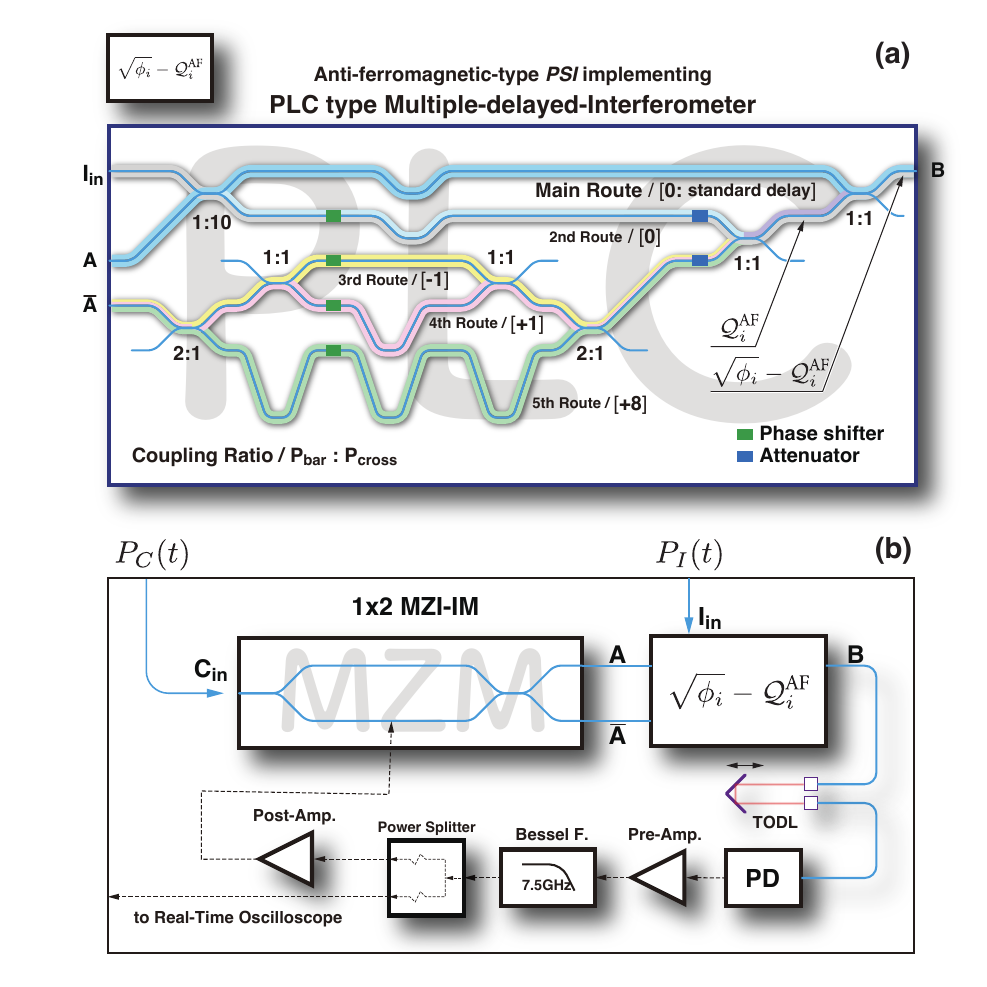}
 \caption{Schematic diagram of the dedicated MaxCut3 ($N\!=\!16$) problem system for the SSBM proof-of-principle experiment. (a) Optical delay interference circuit (ODIC) for physically implementing the pSI in the target problem (see Fig. 3), (b) Configuration diagram of SSBM physically implementing the pSI using the ODIC. The output from the 1$\times$2 MZM's output port ${\bf A}$ and its complementary output port $\overline{\bf A}$ are input to the corresponding input ports of the ODIC. Physical implementation of the pSIs are achieved through delay interference between optical pulses within the ODIC.}
 \end{figure}
Figure 4(a) shows a schematic diagram of the ODIC, and Fig.4(b) shows a schematic diagram of a proof-of-principle system for SSBM that implements this circuit.
Here, the ODIC was fabricated as a silica-based planar lightwave circuit, while the other parts of the system were constructed using commercially available components.
The ODIC has input ports corresponding to the output port ${\bf A}$ and the complementary output port $\overline{\bf A}$ of the 1$\times$2 MZM, and the optical delay interference between the input optical pulses from these input ports physically implements the relationships (pseudo-spin interactions) between the 16 full-dissipative systems, which is then output to the variable optical delay line in the subsequent stage.
This ODIC is designed so that, under the condition that no additional optical amplification is performed, the losses for branching and adjustment (phase, attenuation) are taken into account to satisfy ${\cal J}_{i:k}^{\rm AF}\!=\!1/30\, (k\!=\!i\!+\!8, i\!+\!1, i\!-\!1)$.
The delay of the third route corresponding to $k\!=\!-1$ is set to be the smallest, and the delays of the routes corresponding to the reference route (main route and second route:$k\!=\!0$), $k\!=\!+1$ and $k\!=\!+8$ are set to be $(k\!+\!1)$ times larger than the clock pulse interval.
In addition, by adjusting the phase shifters and attenuators, the pseudo-spin interaction represented by Eqns.(2)-(7) is physically implemented through optical delay interference.
A series of optical clock pulses with a repetition frequency of 1.02\! GHz is extracted from a series of optical output pulses from an active mode-locked laser with a repetition frequency of 8.16 \! GHz and a pulse width of 10 \! ps and then split into two.
To allow repeated trials, optical coherent clock pulses $P_{C}(t)$ with a width of $0.97 \mu s$ and a period of 0.996 \! MHz are extracted from one of these optical pulse trains and reset regions with a width of $33.9 ns$ are set.
From the other of these optical pulse trains, a sequence consisting of the first 16 consecutive optical pulses is extracted and prepared as an optical initializing pulse train.
In addition, we adjusted the timing so that the electrical pulse originating from the first optical pulse in this optical initializing pulse train is input as the intensity-modulated signal to the 1$\times$2 MZM at the same time as the first optical pulse of the optical coherent clock pulses is input to the 1$\times$2 MZM.
With the above settings, SSBM will reset and immediately transition to the NS, causing the SSB phenomenon to start, repeating this behavior.
The peak power of the inflow optical pulse detected by the photodiode corresponds to the transmittance of each system, so this electrical pulse was normalized to monitor the behavior of the system. Using a real-time oscilloscope, a series of 800 trials was recorded at once.
 \begin{figure}[h]
\centering
\includegraphics[width=9cm]{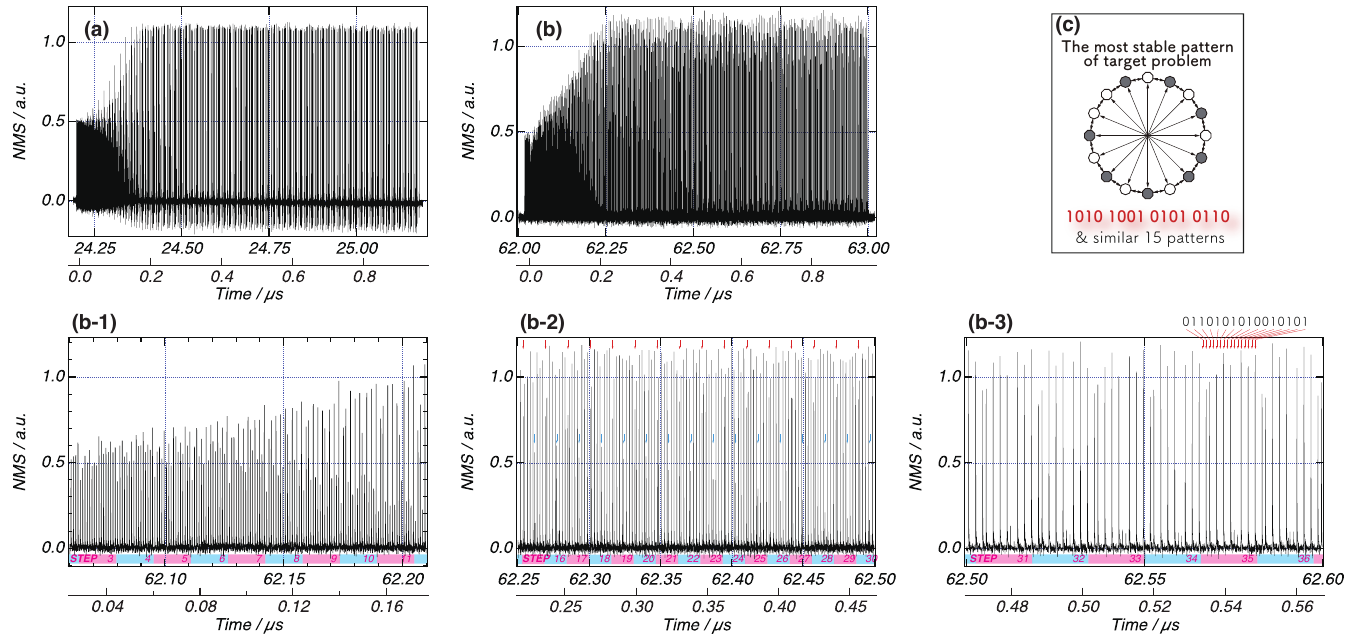}
 \caption{Observed waveforms from the SSBM principle verification experiment. (a): Operational verification waveform with the pSI off (16 independent SSB phenomena observed), (b): Overall observed waveform of one trial with the pSI on, (b-1)-(b-3): Time-domain enlarged waveforms of (b), (c): Example of the most stable state pattern in the target problem. The electrical signal data monitored during 800 repeated trials was acquired in a single batch using a real-time oscilloscope. It can be confirmed that the state created by SSBM corresponds to one of the most stable state patterns of the target problem (see (b-3)). NMS: Normalized Monitoring Signal, , Time axis (lower side): Relative time based on the peak position of the first pulse in the pulse train that started outputting after the reset operation.}
 \end{figure}
The waveforms shown in Fig. 5(a) correspond to the case where the pSI is turned off (set ${\cal J}_{i:k}\!=\!0$) using attenuators in the ODIC.
It was confirmed that the fundamental SSB phenomenon occurs in each of the 16 independent full-dissipative systems, where the system transitions to either the 0-state or 1-state in a manner that reflects the initial fluctuations.
In cases without pSI, since full-dissipative systems are entirely independent and do not compete with each other, state transitions are completed in approximately 20 \! steps, with little variation in the number of steps required.
Figures 5(b)-(c) show the normalized monitoring signals of the SSBM with pseudo-spin interactions.
In this example, it can be confirmed that the state transition is completed in approximately 30 steps, successfully transitioning to the state corresponding to the most stable solution.
Furthermore, the number of steps required for solution search exhibited a relatively large variation and required significantly more steps compared to when pSI was off, with some cases requiring over 50 steps.
This variability is believed to stem from the fact that pSI-induced conflicts exhibit different behaviors in individual cases.
Based on data obtained from experiments, when thresholding was performed defining values above 80\% of the normalized value as {\it up} and values below 20\% as {\it down,} the success rate to find the most stable solution was approximately 97\%. These results demonstrate that SSBM is a promising physical implementation simulator for combinatorial optimization problems. They also strongly support the existence of a duality (correspondence) between the Ising-model-like system states that appear in SSBM and the stable states of the Ising model system.

\section{Numerical Simulation-Based Characteristic Evaluation of SSBM in Large-Scale Problems}
This paper addresses the evaluation of the SSBM's characteristics in larger-scale problems as another important point of issue, and we selected a fully connected complete graph known as $K_{2000}$ \cite{Haribara} as a specific large-scale problem and evaluated it through numerical simulation.
The problem $K_{2000}$ we selected is not only a high-density graph problem where computational complexity increases dramatically with algorithmic solvers, but it also serves as a substantive benchmark problem comparable to prior research reports using other methods \cite{Haribara,Inagaki,Gotoh01,Gotoh02,KYamamoto}.
This section first describes two noteworthy behaviors of the SSBM identified through numerical simulation applied to a relatively small-scale problem ($N\!=\!50$). And, based on the findings obtained, we then propose an improvement method for the SSBM that is physically implementable.

\subsection{Asymmetric Effects in Pseudo-Spin Interactions}
Although the pSI is an interaction defined on the basis of the interpretation of exchange interaction energy, there are clear differences between them. Furthermore, the order parameter itself, which assumes a correspondence with spin, has the distinct difference of being a continuous variable.
Therefore, to gain insight into the resulting correspondence and behavioral differences between the two, we first addressed a relatively small-scale problem ($N\!=\!50$) and performed numerical simulations.
 \begin{figure}[h]
\centering
\includegraphics[width=9cm]{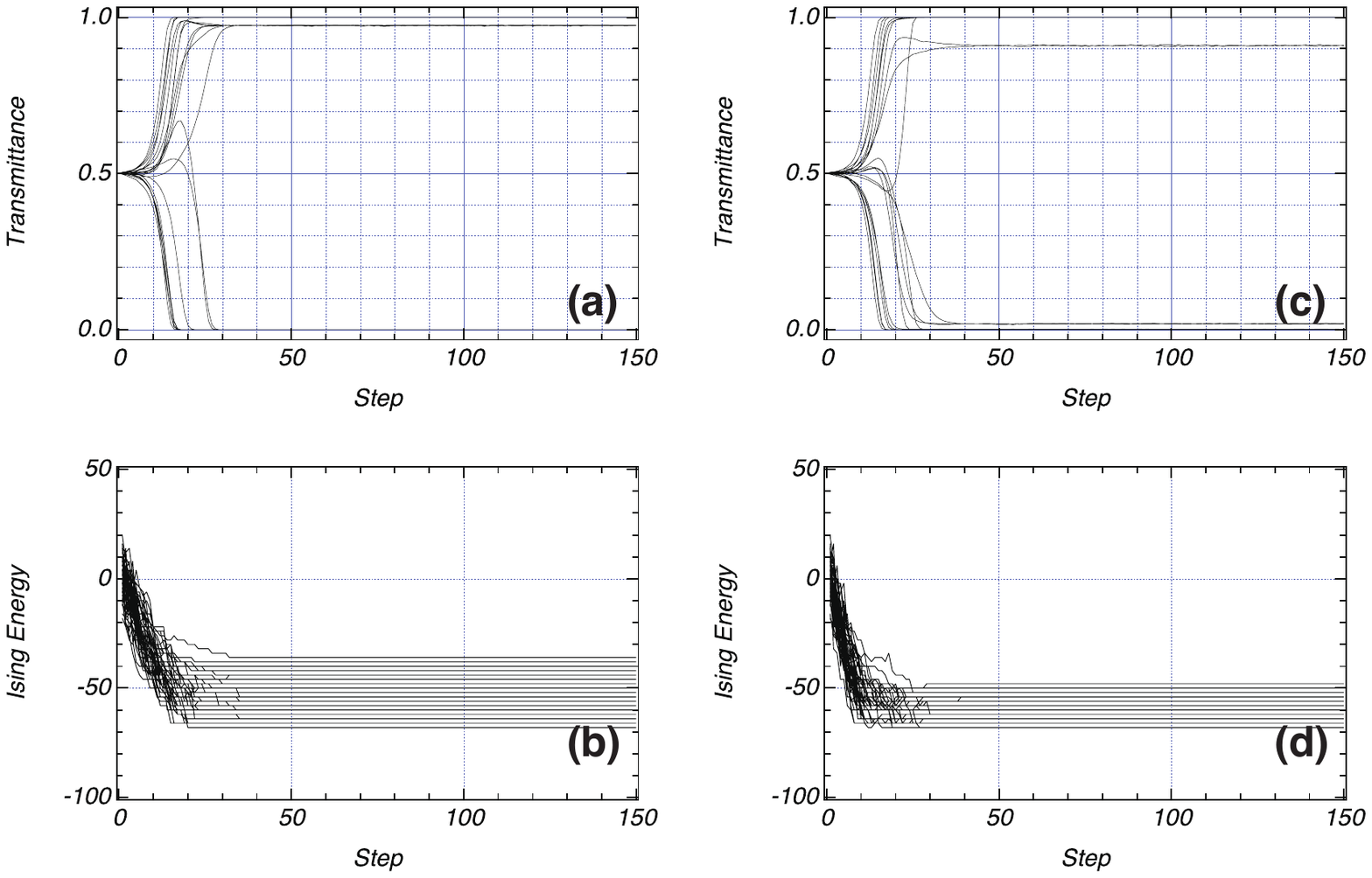}
 \caption{Numerical simulation results of SSBM applied to MaxCut3(N=50). For each of the two cases - ${\cal J}_{i:k}\!=\!0.05$ (case-1) and ${\cal J}_{i:k}\!=\!0.1$ (case-2) - 50 numerical simulations with different initial fluctuations were performed. (a): behavior of the order parameters (pseudo-spins) in case-1 (behavior of 20 pseudo-spins in one simulation sample), (b): behavior of the Ising energies in case-1 (behavior of 50 simulation samples) (c): behavior of order parameters in case-2, (d): behavior of the Ising energies in case-2.}
 \end{figure}
 \begin{figure}[h]
\centering
\includegraphics[width=9cm]{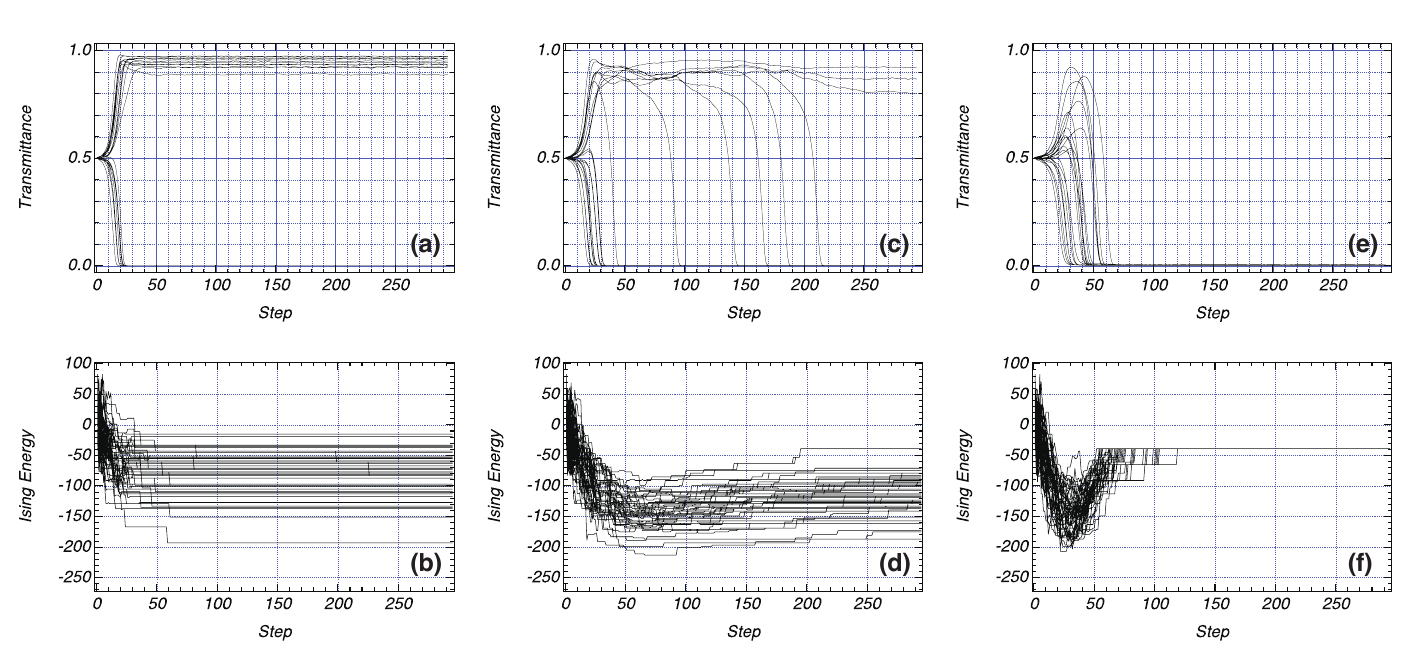}
 \caption{Numerical simulation results of SSBM applied to $K_{50}$. For each of the three cases - ${\cal J}_{i:k}\!=\!0.0025$ (case-1), ${\cal J}_{i:k}\!=\!0.0035$ (case-2), and ${\cal J}_{i:k}\!=\!0.005$ (case-3) -, 50 numerical simulations with different initial fluctuations were performed. (a), (c), (e): Behavior of the order variable (pseudo-spin) in case-1, case-2, and case-3, respectively (behavior of 20 pseudo-spins in a single simulation sample). (b), (d), (f): Behavior of the Ising energy in case-1, case-2, and case-3, respectively (behavior of 50 simulation samples).}
 \end{figure}

Figures 6 and 7 show the results of numerical simulations for the MaxCut3 problem and a fully connected complete graph problem, respectively.
The significant difference in the magnitude of the pSI employed for these two problems is due to the differing densities of coupling, and the calculations were carried out within the variable range of the pSI, focusing on the region where its effects are relatively pronounced \cite{Endnote00}.
As shown in the figure, for the MaxCut3 problem composed solely of anti-Ferro-type interactions, both the behavior of the pseudo-spins (Figs. 6(a), 6(c)) and the behavior of the Ising energy (Figs. 6(b), 6(d)) were confirmed to exhibit the expected transition behavior to the stable state on the Ising energy, consistent with our predictions.
On the other hand, for the complete graph problem ($K_{50}$) where Ferro and Anti-Ferro interactions coexist, it was confirmed that when pSI exceeds a certain threshold (roughly ${\cal J}_{i:k}\!\gtrapprox\!0.0035$), the pseudo-spins are drawn toward the 0-state (Fig. 7(c)), and at ${\cal J}_{i:k}\!=\!0.005$, all pseudo-spins ultimately end up in the 0-state (Fig. 7(e)), and regarding behavior on the Ising energy, while it temporarily stabilizes, it later shifts toward an unstable state and converges to a specific energy level (Fig. 7(d), (f)), exhibiting unexpected behavior.
\begin{figure}[h]
\centering
\includegraphics[width=9cm]{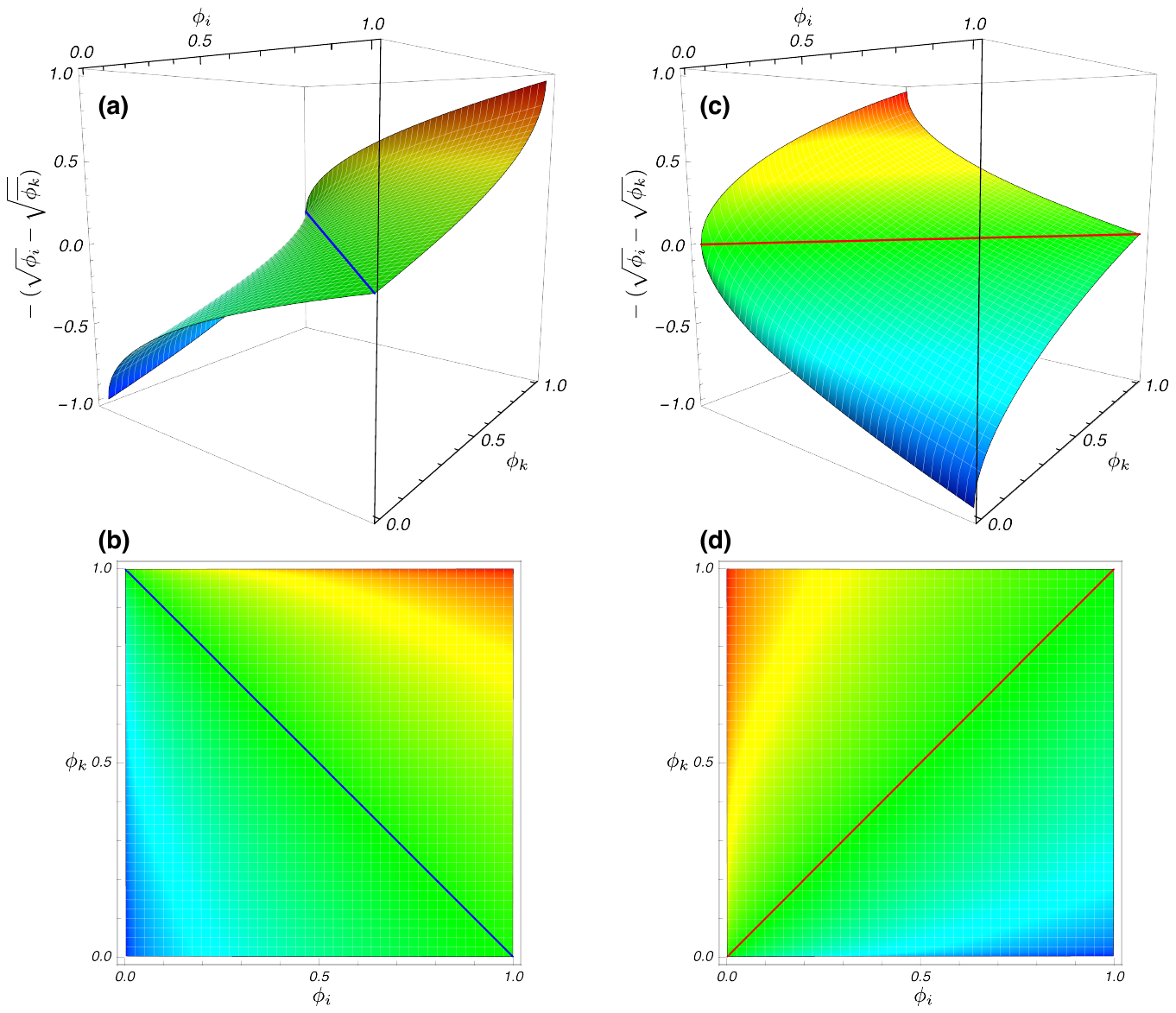}
 \caption{Order parameter (pseudo-spin) dependence of the pSI (between any two pseudo-spins \& ${\cal J}_{i:k}\!=\!1$). (a):Anti-Ferro type, (b):Top view of the graph (a), (c): Ferro type, (d):Top view of the graph (c).}
 \end{figure}
 
Figure 8 graphically illustrates the pseudo-spin ($\phi$) dependence of pSI when ${\cal J}_{i:k} = 1$.
The pSI is defined on the basis of the optical interference phenomenon and is described by the amplitudes of the two interacting pseudo-spins.
Consequently, the pseudo-spin dependence of pSI exhibits the nonlinearity and asymmetry shown in the Fig.8.
The blue and red lines on the graph indicate the positions of $(\phi_{i}, \phi_{k})$ combinations where pSI is 0. Roughly speaking, the effect of pSI on $\phi_{i}$ acts to move from any combination state of $(\phi_{i}, \phi_{k})$ toward the direction of its nearest position of a combination where pSI is 0 (i.e. toward the direction of these blue and red lines).
For Anti-Ferro-type pSI, the pSI acting on $\phi_{i}$ will ultimately be more effective for transitions toward $(\phi_{i}, \phi_{k})$ = (0, 1) than for transitions toward $(\phi_{i}, \phi_{k})$ = (1, 0). This arises because, in the presence of the competing pSI, the asymmetry of the pSI allows the vicinity of $(\phi_{i}, \phi_{k})$ = (0, 1) to approach the Blue-line more rapidly than the vicinity of $(\phi_{i}, \phi_{k})$ = (1, 0) (see the Anti-Ferro-type pSI shown in Fig.8(a) and 8(b)). Consequently, the dynamics-derived attraction becomes more effective.
However, focusing on the transition effect that acts on $\phi_{k}$, the opposite holds true, so this asymmetry cancels out and its effectiveness disappears.
This is considered the reason why the computational simulation results in the MaxCut3 problem exhibit the behavior we initially expected (see Fig.6).
On the other hand, in the case of Ferro-type pSI, the transition effects acting on both $\phi_{i}$ and $\phi_{k}$ are more effective overall when they transition toward $(\phi_{i}, \phi_{k})$ = (0, 0) is more effective than the transition effect toward $(\phi_{i}, \phi_{k})$ = (1, 1) (see Ferro-type pSI shown in Figs. 8(c) and 8(d)), and unlike anti-Ferro-type pSI, there is no mutual cancelation, so the dynamics-derived attraction toward the 0-state becomes more effective.
This can be considered a major factor in the phenomenon in which all pseudo-spins in $K_{50}$ with Ferro-type pSI transition to the 0-state (see Fig.7(e)).
Furthermore, it is thought that the balance between the dynamics-derived attraction toward the 0-state\&1-state and the magnitude of the pSI effect leads to more complex behavior and consequently different results (see Figs.7(a),(c),(e)).
This unexpected phenomenon, while highly intriguing \cite{Endnote03}, it is clear that such effects become a major obstacle when considering the application of SSBM to combinatorial optimization problems.
Moreover, if the above predictions and hypotheses are correct, then a solution can be proposed relatively easily.
Specifically, it is thought that this asymmetry can be resolved by constructing a system conjugate to the one examined thus far and averaging it with the original system. The crucial point here is that both this conjugate system and the system averaged together with the conjugate system are physically implementable, just like the conventional SSBM.
Here, the equations governing this composite system can be described as follows.
\begin{eqnarray}% equation no. 08
\phi_{i+m}
& = &  {\cal C}( \,\frac{\pi}{2}, \,\phi_{i},  \,{\cal Q}_{i} \,)\\
& = &
\left\{\sin^{2}(\frac{\pi}{2} \,|\sqrt{\phi_{i}} -  {\cal Q}_{i} |^{2} ) +  \overline{ \sin^{2}(\frac{\pi}{2} \,|\sqrt{\overline{\phi_{i}}} -  {\overline{\cal Q}_{i} } |^{2} )} \right\}/ \,2 \notag 
\end{eqnarray}

\begin{figure}[h]
\centering
\includegraphics[width=9cm]{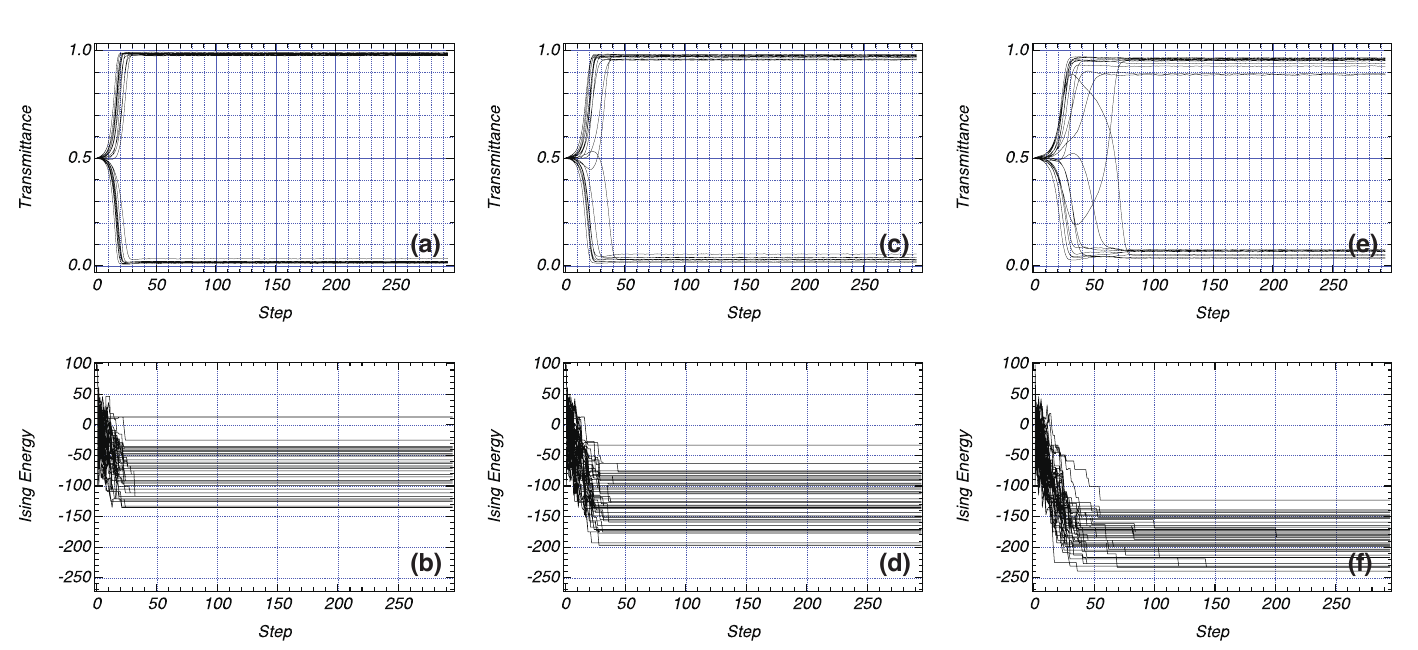}
 \caption{Numerical simulation results of the improved SSBM (SSBM following Eq.(8)) applied to $K_{50}$. For each of the three cases - ${\cal J}_{i:k}\!=\!0.0025$ (case-1), ${\cal J}_{i:k}\!=\!0.0035$ (case-2), and ${\cal J}_{i:k}\!=\!0.005$ (case-3) -, 50 numerical simulations with different initial fluctuations were performed. (a), (c), (e): Behavior of the order variable (pseudo-spin) in case-1, case-2, and case-3, respectively (behavior of 20 pseudo-spins in a single simulation sample). (b), (d), (f): Behavior of the Ising energy in case-1, case-2, and case-3, respectively (behavior of 50 simulation samples).}
 \end{figure}

Under the same conditions (see Figs.6 and 7), the verification of numerical simulation based on Eq.(8) revealed that while the case applied to MaxCut3 showed almost no difference from the conventional method, the case applied to $K_{50}$ exhibited a drastic difference.
Figure 9 shows the results of the numerical simulations concerning $K_{50}$. The phenomenon of avalanche-like transition to the 0-state observed in Fig.7(c)\&(e) has disappeared, and stabilization, as measured by the Ising energy, was shown to proceed without hindering.
These results clearly demonstrate that the aforementioned predictions and hypotheses are correct.
Furthermore, in the following subsection, we will continue our discussion based on the composite system described by Eq.(8).

\subsection{Competition between dynamics-derived attractions and pSI}
\begin{figure}[h]
\centering
\includegraphics[width=9cm]{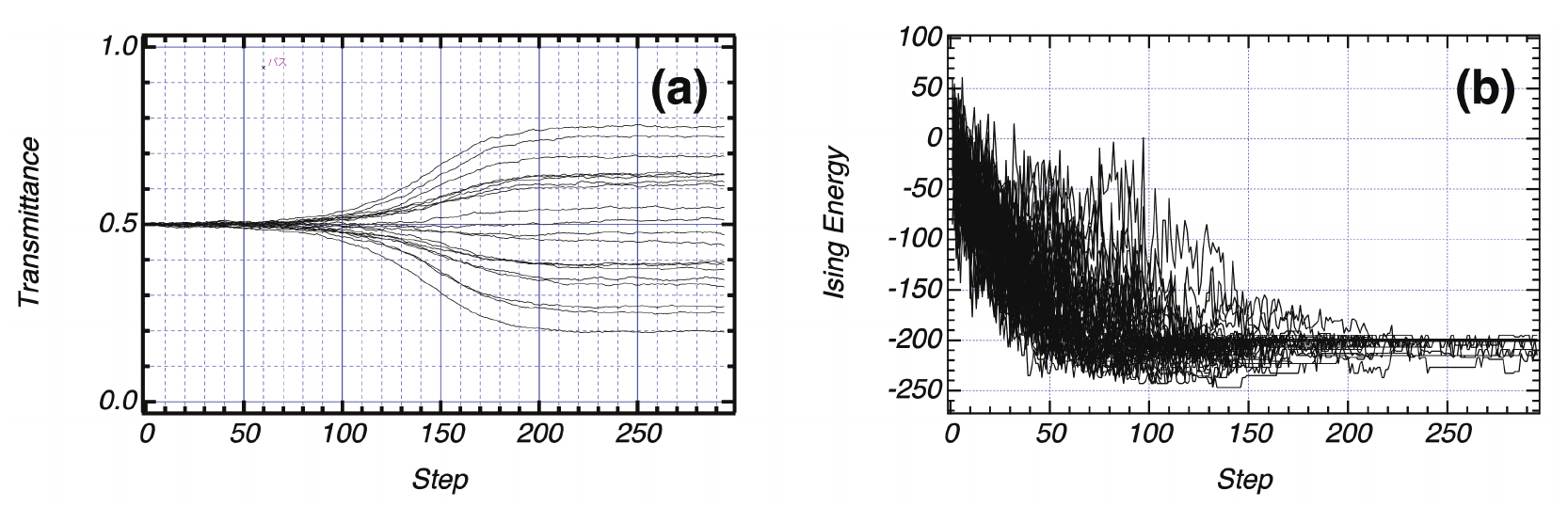}
 \caption{Numerical simulation results of SSBM following Eq.(8) applied to the $K_{50}$. We performed 50 numerical simulations with different initial fluctuations for a case where the pSI was set significantly larger than in the case shown in Figure 9 (${\cal J}_{i:k}\!=\!0.0095$). (a):Behavior of the order parameters (pseudo-spin) (behavior of 20 pseudo-spins in a single simulation sample). (b):Behavior of the Ising energy (behavior of 50 simulation samples).}
 \end{figure}

Figure 10 shows the numerical simulation results for $K_{50}$ obtained using Eq.(8), illustrating an example where the pSI was set quite large.
Note that the behavior of the Ising energy with respect to the step (see Fig.10(b)) was calculated using results where $\phi \!\geqq\! 0.5$ was classified as 1-state and $\phi \!<\! 0.5$ as 0-state.
As can be seen from these figures, the competition between the attraction derived from dynamics and pSI causes the pseudo-spin to fail to sufficiently attract toward the 0-state/1-state in regions where the pSI value is large.
Unlike Ising model-based simulators (such as coherent Ising machine or simulated bifurcation machine), SSBM guarantees behavior based on the interpretation of exchange interaction energy even in such situations. Therefore, it is possible to promote stabilization as measured by the Ising energy, and the behavior shown in Fig.10 also demonstrates this fact.
However, when considering the practical operation of SSBM, it is anticipated that each pseudo-spin will become an obstacle in the process of determining whether it converged to the 0-state or 1-state.
Furthermore, the fact that the stable state where the pseudo-spin is extremely far from the 0-state \& 1-state is displaced from the Ising energy potential surface raises concerns about the degradation of performance as a simulator for the COP, which requires finding stable states as Ising energy.
Therefore, we investigated ways to enhance the attraction to 0-state\&1-state arising from dynamics and evolve the system into one capable of achieving proper state separation of pseudo-spins (namely, ensuring that the pseudo-spin converges sufficiently to the 0-state\&1-state) against competing pSI.
(Although studies on improving convergence to up-state and down-state have been reported for Ising machines utilizing bifurcation phenomena \cite{Bohm02, Shi}, the SSBM, which employs dynamics-derived attraction, required an entirely different approach.)
As a result, we arrived at an evolved version of SSBM expressed by the following equations.
\begin{eqnarray}% equation no. 09 ~ 10
\phi_{i+m}
& = &
{\cal N}(\,n, \,\pi, \, {\cal C}( \,\frac{\pi}{2}, \,\phi_{i},  \,{\cal Q}_{i} \,))\\
{\cal N}(n, \pi, P_{in} )
& := &
\underbrace{
\sin^{2}(\frac{\pi}{2}\,\sin^{2}(... \frac{\pi}{2}\,\sin^{2}(}_{n:\rm{{ number\, of\, nesting}}}\frac{\pi}{2}\,P_{in} )...))
\end{eqnarray}

It should be noted that the evolved SSBM proposed here also satisfies the conditions of being physically implementable and observable as a phenomenon.
\begin{figure}[h]
\centering
\includegraphics[width=9cm]{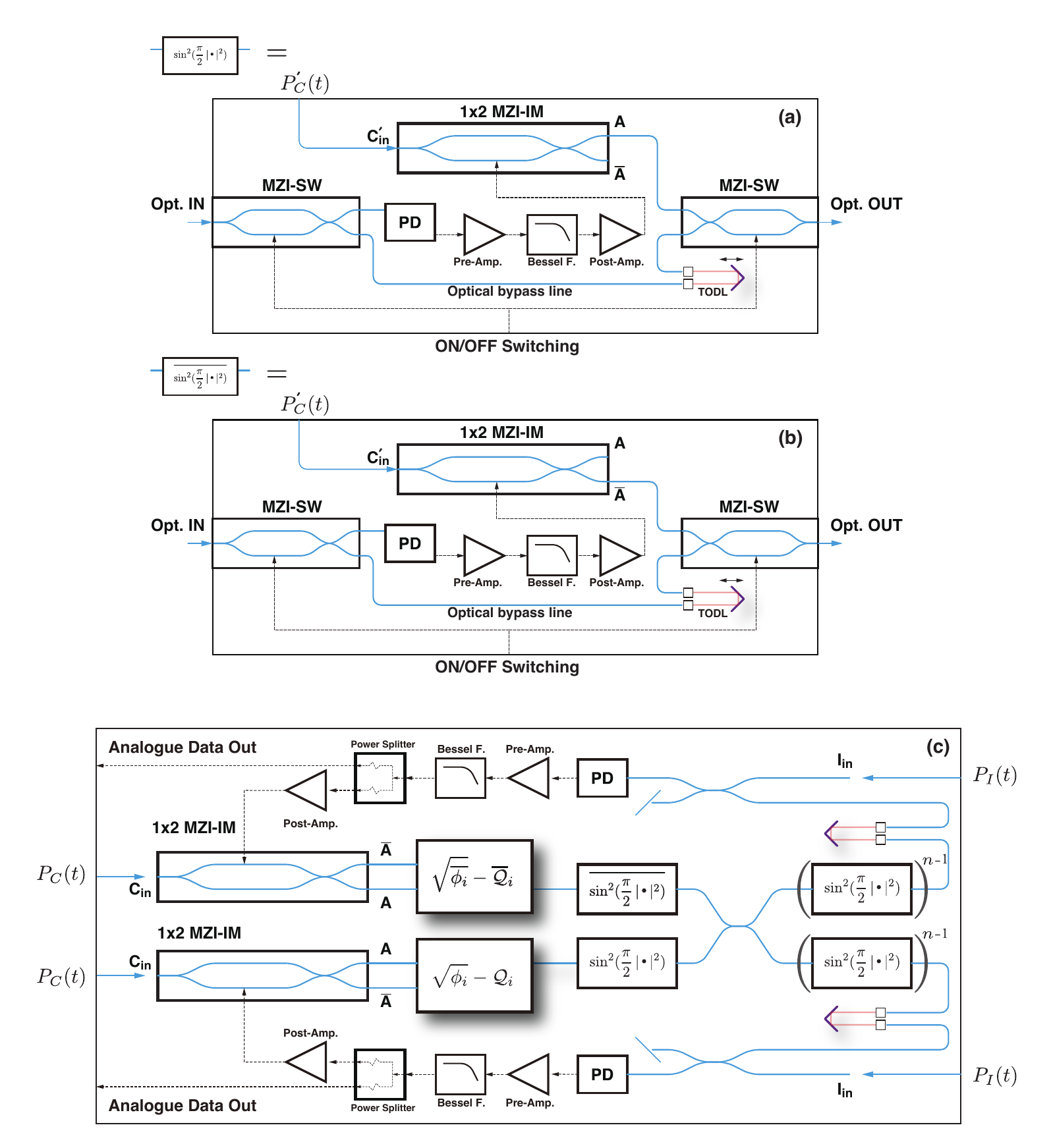}
 \caption{Schematic diagram of the circuit configuration for the evolved SSBM (SSBM following Eq.(9)). (a): Circuit unit used to realize the nested action of the iteration equation governing the system. (b): Circuit unit outputting the complement of the nested action. These two circuit units are supplied with optically synchronized clock pulses. It also includes a mechanism to switch to a bypass optical waveguide adjusted for the same delay, enabling dynamic ON/OFF switching of this effect. (c): Overall circuit configuration of the evolved SSBM. Two optical clock pulse trains and two initialization optical pulse trains are supplied with their timing adjusted.}
 \end{figure}
Specifically, this evolved SSBM can be realized using the optical system with the circuit configuration shown in Fig.11.
The circuit unit shown in Fig.11(a) is used to embody the nested action of the iterative equation followed by the system,
and the number of nested actions (NNA) $n$ can be increased by increasing the number of circuit units connected in a cascade (Fig.11(c)).
And, this circuit unit is designed to enable switching between the functional circuit section and the optical bypass line, anticipating dynamic changes to the NNA.
\begin{figure}[h]
\centering
\includegraphics[width=9cm]{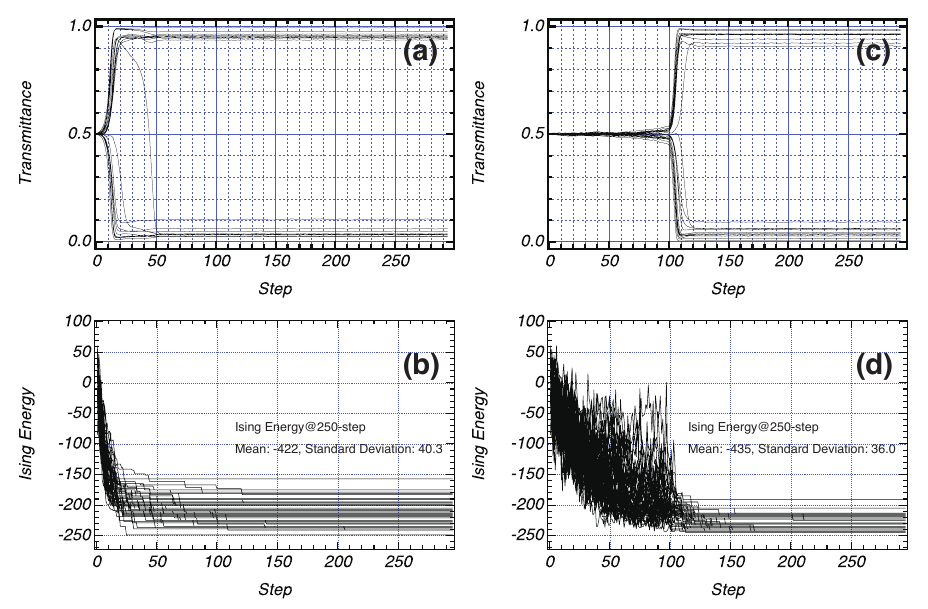}
 \caption{Numerical simulation results for the evolved SSBM (SSBM following Eq.(9)) applied to $K_{50}$. Cases where number of nested actions ($n$) was set to $n\!=\!2$ ((a) and (c)) and where $n\!=\!1$ for the first 100\! steps and $n\!=\!2$ thereafter ((b) and (d)) were examined. Fifty numerical simulations with different initial fluctuations were performed with ${\cal J}_{i:k}\!=\!0.0095$. (a) and (c): Behavior of the order variable (pseudo-spin) (behavior of 20\! pseudo-spins in a single simulation sample). (b) and (d): Behavior of the Ising energy (behavior of 50\! simulation samples).}
 \end{figure}
Figure 12 shows the results of the numerical simulation for the case where the NNA was controlled.
In the case where the number of nested actions was increased by one (fixed), it was confirmed, as initially intended, that the separation of the 0-state \& 1-states of the pseudo-spin could be effectively achieved (Fig.12(a)) and that stabilization in terms of Ising energy was improved (Fig.12(b)).
Furthermore, in this case, it was observed that the variance in the Ising energy in the final state (step$\sim$300) was relatively large.
This is thought to occur because increasing the NNA rapidly advances both the convergence of the pseudo-spins toward the attractor and the fixation of the state.
On the other hand, in the case of dynamic control where the NNA is not increased until step=100 and increases once thereafter, not only the state separation of pseudo-spins in the final state (Fig.12(c)) and the improvement in stabilization (Fig.12(d)), but also the suppression of variation in Ising energy was confirmed (Fig.12(d)).
This can be considered evidence that, during the initial 100 steps, there is not only the disadvantage that pseudo-spin state separation does not progress, but also the advantage that transitions to more stable pseudo-spin configurations are occurring gradually yet reliably.
Furthermore, it is thought that this advantage is realized by pSI-specific effects not present in exchange interactions, suggesting that effective control of NNA could enable more advanced stabilization control (solution exploration).

\subsection{Numerical Simulation of SSBM Applied to $K_{2000}$ Problem}

Here, we describe numerical simulations applying the aforementioned evolved SSBM to the $K_{2000}$ problem \cite{Haribara}, which has been used as a benchmark problem for Ising machines \cite{Inagaki,Gotoh01,Gotoh02,KYamamoto}.
Since $K_{2000}$ is a fully connected graph, the upper bound of the value of ${\cal J}_{i:k}$ is approximately determined by the reciprocal of the number of pSIs associated with the element $\phi_{i}$ (see Eq.(4), Eq.(5)).
And, to maximize the pSI effect, we set ${\cal J}_{i:k}$ close to this upper bound and adjusted the NNA and number of steps to ensure SSBM functions under this condition.

\begin{figure}[h]
\centering
\includegraphics[width=9cm]{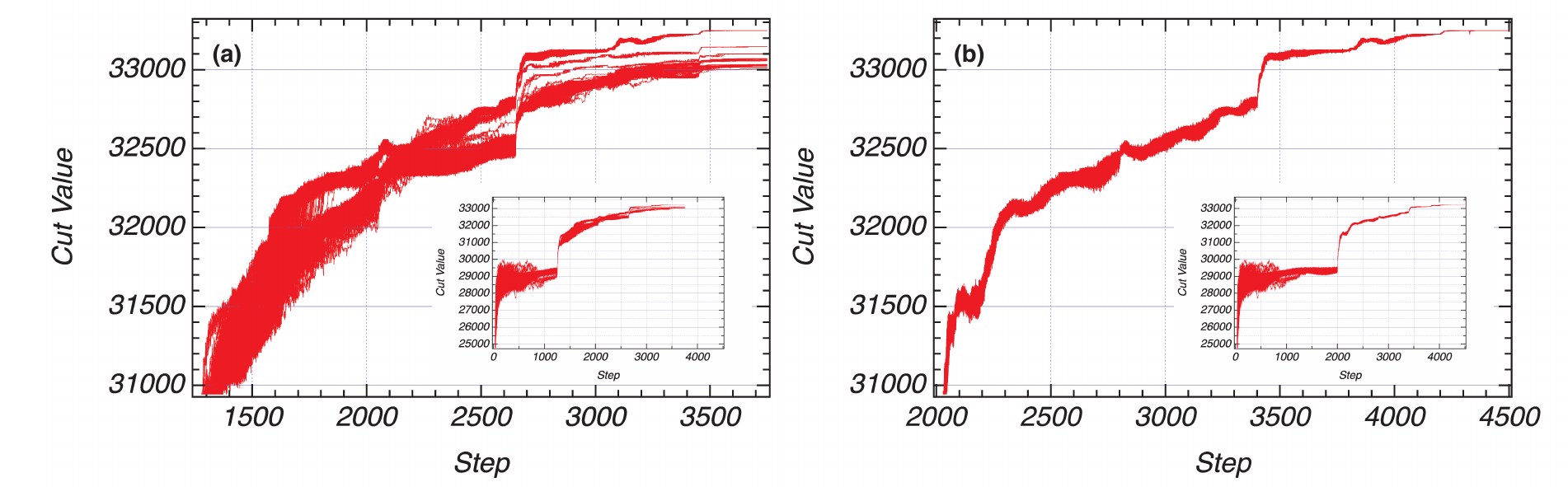}
 \caption{Numerical simulation results (cut valuses behavior) of the evolved SSBM (SSBM following Eq.(9)) applied to $K_{2000}$. Common conditions:Number of simulation samples 1000, ${\cal J}_{i:k}\!=\!4.96\pm0.03\!\times\!10^{-4}$, number of nested actions incrementally increased from $n\!=\!6$ to $n\!=\!12$. (a):(case-1) Set the number of steps to 1250, 800, 600, 400, 400, and 300 for $n\!=\!6,7,8,9,10$, and $12$. (b):(case-2) Set the number of steps to 2000, 800, 600, 400, 400, and 300 for $n\!=\!6,7,8,9,10$, and $12$.}
 \end{figure}
 
 Figure 13 shows the change in cut values with a step progression for 1000 numerical simulations. Since many reports dealing with large-scale problems use cut values rather than Ising energys, we followed this convention.
Additionally, as calculation conditions, the NNA was incrementally increased from $n\!=\!6 $ to $n\!=\!12$, and the value of ${\cal J}_{i:k}$ was adjusted to gradually decrease within a narrow range of $4.96\pm0.03\times10^{-4}$ for each NNA value, and Gaussian random numbers \cite{Matsumoto} (with a standard deviation of $1.0\times10^{-6}$) differing for each simulation were used as fluctuations.
Specifically, in the case of Fig.13(a) (case-1), for $n\!=\!6, 7, 8, 9, 10$ and $12$, the number of steps was set to $1250, 800, 600, 400, 400$ and $300$, respectively. Then in the case of Fig.13(b) (case-2), only the number of steps for $n\!=\!6$ was increased from $1250$ to $2000$.
And, Figs.14(a)-(c) show the histograms of the cut values for 1250, 1400, and 3750 steps in case-1, respectively, while Figs. 14(d)-(f) show those for 2000, 2150, and 4500 steps in case-2, respectively.
The simulation parameters used here were identified through repeated trials. Future challenges include developing methods to find parameters more efficiently using AI technology or other techniques, as well as examining their applicability to a broader range of problem types.
These simulation results demonstrate that SSBM possesses distinct characteristics that distinguish it from the other COP solvers.
For example, in the case of quantum computers (quantum annealers), the tunneling effect can only be expected with a certain probability within a realistic computational time (finite time). Consequently, statistical variation in the obtained solution (cut values) is unavoidable. Furthermore, statistical variation in the cut values obtained has also been confirmed in reported cases involving other Ising machines, including algorithmic types \cite{Quinton,Inagaki,Gotoh02,KYamamoto}.
On the other hand, SSBM exhibited an extremely unique and interesting behavior in which a single stable state is created (Figure 13(b)), and the cut value observed there indicated that SSBM has top-level characteristics (99.7\% of the best value currently known \cite{Gotoh02}). And, this single stable state was confirmed to be paired with a state where the 1-state and 0-state are exactly reversed, and it was also confirmed that both states could be obtained.

We consider this distinctive characteristic of SSBM to be significant in two respects, one of which is as a key to deepening our understanding of the phenomenon itself (its scientific value).
We believe that the SSB phenomenon (process) employed in SSBM, and furthermore the newly introduced pSI in SSBM, create a significant difference from the Ising machines.
For example, in Ising machines, where continuous variables correspond to spins, gradually varying parameters allow the system to progress toward a stable state while remaining close to equilibrium \cite{Inagaki,Gotoh01,Gotoh02}. Broadly speaking, this operation separates all the continuous variables into two groups, up and down, effectively maintaining conditions that keep the system near the Ising energy potential surface. This operation is necessary to induce the transition to a stable state in systems where exchange interactions are introduced directly.

On the other hand, in SSBM, the pSI can function effectively regardless of the relative positional relationship of continuous variables (even when they are not separated into two distinct states corresponding to up and down).
 And, the attraction toward the 1-state or 0-state inherent in the dynamics themselves competes with pSI, determining the stable state as a pseudo-spin. Consequently, for example, at $n\!=\!6$, it exhibits extremely gradual changes or stagnation while assuming a scattered state over a wide range between the 1-state and 0-state.
Moreover, the gradual changes observed here suggest a transition toward a stable state in a particular sense, even though they do little to further advance stabilization from the perspective of the Ising energy.
\begin{figure}[h]
\centering
\includegraphics[width=9cm]{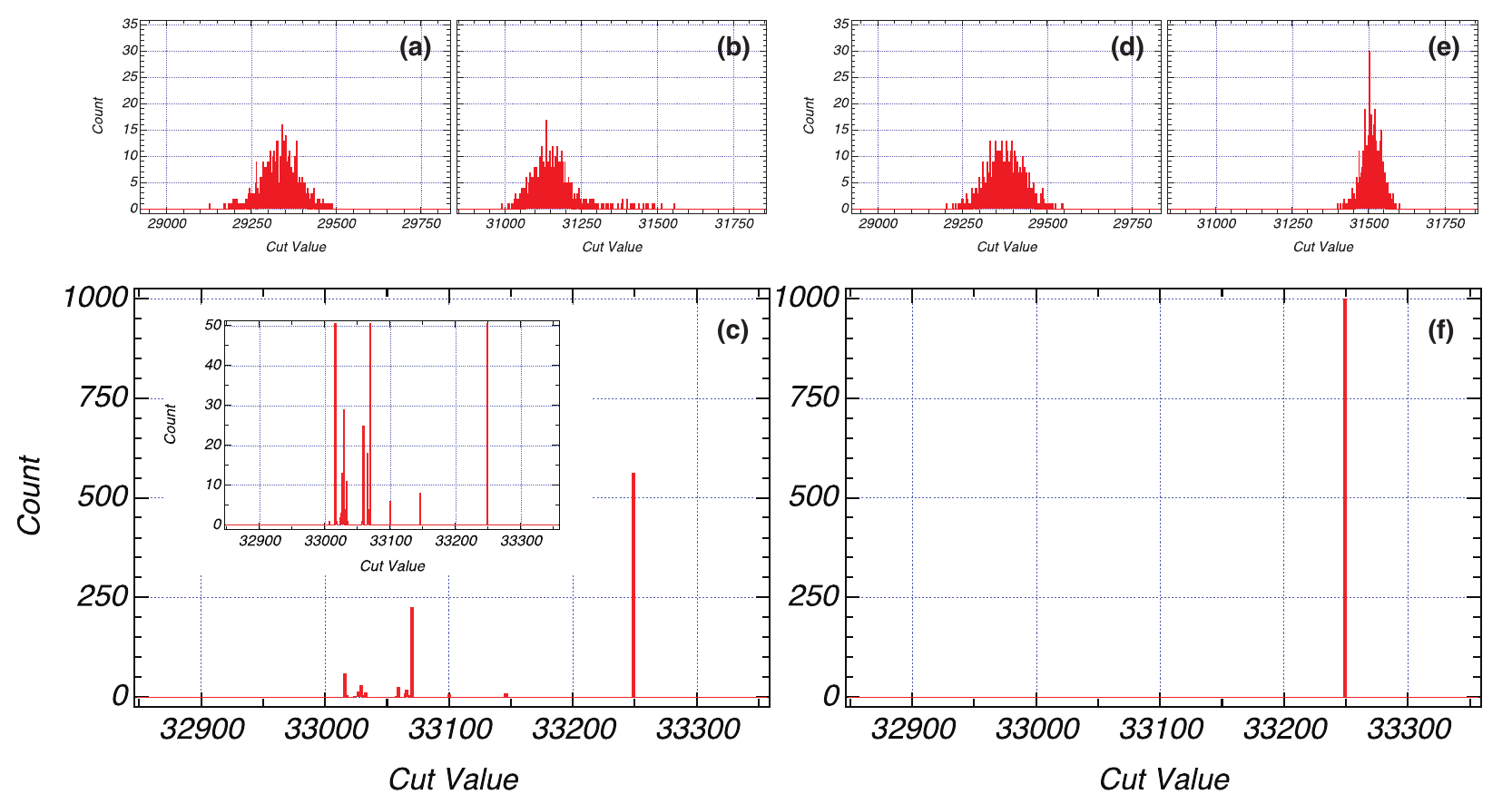}
 \caption{Numerical simulation results (histogram of the cut values) for the evolved SSBM (SSBM following Eq.(9)) applied to $K_{2000}$. (a)-(c):Numerical simulation results under the same conditions as Fig.13(a), showing histograms of cut values at steps 1250, 1400, and 3750, respectively. (d)-(f):Numerical simulation results under the same conditions as Fig.13(b), showing histograms of cut values at steps 2000, 2150, and 4500, respectively.}
 \end{figure}
For example, while there is little difference in the cut values at the final step for $n\!=\!6$ between case-1 and case-2 (cf. Fig.14(a) and Fig.14(d)), a significant difference becomes apparent when comparing the 1400-th step of case-1 (Fig.14(b)) with the 2150-th step of case-2 (Fig.14(e)) (150 steps after switching the NNA to $n\!=\!7$).
In the final step, it was confirmed that the stable state obtained in case-1 retains statistical fluctuations (Fig.14(c)), while convergence to a single stable state was observed in case-2 (Fig.14(f)).
Furthermore, the following points were also confirmed from the simulation results.

{\bf Result-1:}  
Regarding the sample group for case-2, in the stage shown in Fig.14(e), the peaks are sharper compared to the stage in Fig.14(d), but no convergence occurred as a combination pattern, with 98.9\% being mutually different patterns.

{\bf Result-2:}  
Regarding the sample group for case-1, in the stage shown in Fig.14(b), it is nearly separated from case-2 (cf. Fig.14(e)), but in the final step, over 50\% converges to the same state as the single stable state in case-2.

{\bf Result-3:}  
The behavior of cut value (Ising energy) (Fig.13(a)) reveals that there are streams converging to several specific stable states and that the transition is not merely monotonically toward states with larger cut values (states with lower Ising energy).

{\bf Result-4:}  
As NNA increases, continuous variables progressively separate into two distinct states (up and down) with greater clarity, while the combination patterns gradually consolidate.

Taking into account the consistency of the results confirmed above, the following is anticipated.

{\bf Prediction-1:} (Based on {\bf Result-1}  and {\bf Result-2})
The combination patterns leading to a single stable state are distributed considerably more widely than the distribution width observed in Fig.14(e), and their number is significantly larger than that of combination patterns leading to other stable states.

{\bf Prediction-2:}
The combination patterns forming the streams described in {\bf Result-3} broadly constitute a shape resembling a river with many tributaries (or a tree structure), composed of a sequence of combination patterns that differ slightly from each other. The transition effect realized by the pSI at $n\!=\!6$ is effective for traversing between these streams toward a stream ({\bf stream-one}) leading to a single stable state.

{\bf Prediction-3:}
(Consideration that the cut value histogram in the 1450-th step for n=6 in case-2 is identical to Fig.14(a))
Although no significant change in cut value (Ising energy) is observed during the transition from Fig.14(a) to Fig.14(d), it is considered that the state change shifting to {\bf stream-one} is completed during this period. This state change is driven by the entropy effect arising from the relatively large number of states constituting {\bf stream-one} within this energy range.

In an analog variable-type Ising machine that directly inherits Ising interactions using branching phenomena, operations are typically performed to slowly increase the gain and feedback strength so that the steady-state approximation or adiabatic approximation holds. This maintains a state separated into two states, upper-side and lower-side, while approaching the up-spin or down-spin state. In such a process, the system does not fall into the state shown in Figure 10(a). However, if forced into a similar configuration, the interaction itself exhibits behavior resembling significant noise. The system then undergoes a behavior that resolves this state and settles into the separated upper-side and lower-side states. In contrast, SSBM, by introducing a pseudo-spin interaction, guarantees behavior consistent with the interpretation of exchange interaction energy even when the system is not separated into two states (e.g., as in Figure 10(a)), maintaining this state. Noting that the only difference in conditions is the number of steps for n=6, the above results suggest that the unique characteristic of SSBM in the region not separated into these two states effectively led to pattern convergence by ensuring a sufficient number of steps.
Verification from a different perspective to reinforce the above predictions remains a future task.

Another important consideration is the advantages and disadvantages of its application as a COP solver.
For example, when considering application scenarios for the COP solvers (quantum anealer, coherent Ising machine, simulated bifurcation machine and STATICA etc.) where the obtained solutions exhibit statistical variability, it is necessary to determine the number of computational samples and further consider the need for subsequent comparison and selection processes for the obtained sample solutions.
And, it will also be important to determine the extent to which the device configuration can support parallel computing.
On the other hand, since SSBM holds promise for effectively converging patterns that can be obtained under certain conditions, we anticipate that gaining broader insights into the relationship between simulation conditions and convergence effects could lead to discovering efficient utilization methods that conserve computational resources. While it is true that the solutions obtained via SSBM have not yet reached the current best cut value, we consider it a promising option since the solutions achieved represent top-level performance (99.7\% of the currently known best value).
Furthermore, it is important to note that the scientific understanding of the phenomena leading to this single solution and the examination of its impact on applications remain major challenges for the future.

\section{Conclusion}

We have successfully demonstrated experimentally the operation of SSBM and confirmed its promise as a new physically implemented type COP solver.
 This result also strongly supports the existence of a duality (correspondence) between the Ising-model-like system states that appear in SSBM and the stable states of the Ising model system.
In addition, the characteristics were evaluated by numerical simulation on a benchmark problem ($K_{2000}$), assuming application to large-scale problems.
The results of 1000 simulations with different fluctuations confirm that SSBM is able to search for states corresponding to 99.7\% of the best cut value currently known, without statistical variability.
This is a unique characteristic rooted in SSBM's principles, not found in other solvers, and has the potential to become a distinct advantage in the future.
However, in configurations that realize all individual pSIs presented here via optical interference, as the number of pSIs increases, the optical power of each individual pSI must be reduced, making the SN ratio of the optical pulses the root cause of a severe scalability problem. Similarly, for coherent Ising machine, a solver also using optical pulses, no experimental reports exist of configurations achieving Ising interactions solely through optical interference, even for problems where all spins are coupled. This is because coherent Ising machine faces challenges similar to those of SSBM. This scalability problem is expected to be solvable by following the example of coherent Ising machine, which succeeded in experimental application to the fully-connected problem with N=100000 \cite{Honjo}. This scalability issue is expected to be resolved by following the example of coherent Ising machine, which successfully applied to the N=100,000 fully connected problem. There, the problem is addressed using a method termed a measurement and feedback scheme: optically encoded signals corresponding to spins are converted to electrical signals, then interactions are computed digitally using an FPGA, and finally the computed interactions are feedback as optical signals via electro-optical conversion.

\section{Acknowledgment}
The author thanks A. Sugiyama,  M. Tsuda, S. Oka, H. Takenouchi, T. Komatsu, T. Sawada, O. Kagami, S. Mutoh, T. Haga, and A. Okada for their support.


\begin{thebibliography}{43}
\bibitem{Endnote01} 
The virtual boundary referred to here denotes something that cannot be defined solely by spatial conditions. It is set by focusing on the pair of pulses—the $i$-th pulse and the ($i\!+\!m$)-th pulse—among the supplied clock pulses, where the earlier pulse influences the later pulse via the MZM. In other words, it is a boundary that includes temporal conditions, not merely stating that they share the same spatial path. In the system under consideration here, $m$ independent virtual boundaries are set within the same space.
\bibitem{Endnote04} 
The phenomenon of interest in this study is that the system cannot be described by the iteration relation shown in Equation (1) unless it satisfies the special conditions described in Section 2. This is because, unlike systems based primarily on digital processing, the system's information gradually dissipates due to the accumulation of effects such as fluctuations, no matter how precisely the initial conditions are adjusted. Conversely, satisfying these special conditions establishes the relationship shown in Equation (1), namely, robust causal relationships between steps. This paper refers to this causal relationship as dissipative causality.
\bibitem{Endnote02} 
The pseudo-spin interaction in this COP can be realized without dynamic switching, as shown in the optical delay interferometer circuit in Fig. 4(a), due to the regularity of the interacting elements. However, to achieve interaction when no regularity exists, it becomes necessary to dynamically turn ON/OFF pre-prepared optical delay interferometers, making both the circuit itself and its control more complex.
\bibitem{Endnote00} 
For example, in MaxCut3, since “each element interacts with the other three elements,” the pSI must be kept at approximately one-tenth or less. This is because if the pSI becomes too large, the process by which the pseudo-spin converges to the attractor is significantly hindered, preventing the attainment of a stable state.
\bibitem{Endnote03} 
This result predicted by numerical simulations — namely, that in a many-body-like system created by symmetry breaking, dual elements expected to be nearly balanced in population instead exhibit behavior where only one element is created — is not only reminiscent of the physical process called baryogenesis (an unsolved theoretical research subject in physical cosmology), but is also particularly intriguing as it represents an observable singular case.
\end{thebibliography}
\end{document}